\documentclass[11pt,a4paper,leqno]{article}
\usepackage{latexsym}
\usepackage{amsfonts}
\usepackage {amsbsy}
\usepackage {amsmath}
\usepackage{amssymb}
\newtheorem{defn}{Definition}[section] 
\newtheorem{thm}[defn]{Theorem} 
 
\newtheorem{lem}[defn]{Lemma} 
\newtheorem{cor}[defn]{Corollary}

\newcommand \psh {plurisubharmonic }
\newcommand \demo { Proof: }
\newcommand \C{\mathbb C}

\newcommand \R{\mathbb R}
\newcommand \B {\mathbb B}
\newcommand  \G {\mathbb G}

\newcommand \mF {\mathcal F}

\newcommand \mL {\mathcal L}
\newcommand \mP {\mathcal P}
\newcommand \mU {\mathcal U}
\newcommand \mcal {\mathcal}
\newcommand \fin{$\blacktriangleright $\\}
\newcommand \vep  {\varepsilon}
\newcommand \m  {\mathcal}
\newcommand \mb {\mathbb}
\newcommand \W {\Omega}

\newcommand \la {\lambda}
\newcommand \vphi {\varphi}
\newcommand \Sub {\Subset}
\newcommand \sub {\subset}
\newcommand \sm {\setminus}
\newcommand \ove {\overline}
\newcommand \lgra {\longrightarrow}
\newcommand \mrm {\mathrm}
\newcommand \es {\emptyset}
\numberwithin{equation}{section}
 
 \title{Polya's inequalities, global 
 uniform integrability and the size of \psh lemniscates}
 \author{ S. Benelkourchi, B. Jennane
  and A. Zeriahi
 \noindent \footnote {This work was partially supported by
 the programmes PARS MI 07 and AI.MA 180}} 
\date{}
\begin{document}
\maketitle

\noindent {\bf Abstract} First we prove a new
 inequality comparing uniformly the
 relative volume of a Borel subset with respect to
 any given  complex euclidean ball $\B \sub \C^n$
 with its relative logarithmic capacity in $\C^n$ with
 respect to the same ball $\B$.
 An analoguous comparison inequality for Borel subsets
 of euclidean balls of any generic real subspace of $\C^n$
 is also proved. 
 
 Then we give  several interesting applications of 
 these inequalities. 
 First we obtain sharp uniform estimates on the relative size 
 of \psh lemniscates associated to the Lelong class 
 of \psh functions of logarithmic singularities at 
 infinity on $\C^n$  as well as  the Cegrell class of 
 \psh functions of bounded Monge-Amp\`ere mass 
 on a hyperconvex domain $\W \Sub \C^n.$ 
 
 Then we also deduce new results
 on the global behaviour of both the Lelong class 
 and the Cegrell class of \psh functions.

\section{Introduction}
 Local uniform integrability and estimates on the 
 size of sublevel sets of \psh functions in terms of 
 capacities or various measures have been studied 
 earlier in several works (cf. [Cu-Dr-Lu], [Ki], [Ko 2], 
 [Ze 2], [Ze 3], [Pl], [Be-Je]). 
  Such estimates turn out
 to be useful in many areas of Complex Analysis 
 as Pluripotential Theory,
 Pad\'e Approximation and Complex Dynamics 
 (cf. [Ki], [Ko 1], [Ko 2], [Cu-Dr-Lu], [Fa-Gu]). 

 Our aim here is to generalize the classical Polya's 
 inequality to subsets of any generic subspace 
 of $\C^n$ and to give several new applications to 
 the study of the global behaviour of two 
 important classes of \psh functions.

 More precisely, given a generic subspace
 $\mb G \sub \C^n,$ we prove a new inequality 
 estimating from above the
 relative volume in $\G$ of a Borel subset 
 with respect to an euclidean ball 
 $B \sub \mb G $  in terms of its relative 
 logarithmic capacity in $\C^ n$ with
 respect to the same ball $B$, up to a multiplicative 
 numerical constant which depends only on the dimension of
 $\G$ but not on the "condenser" considered.
   
 Formulated in this way, Polya's inequalities turn 
 out to play an important role in applications, 
 implying interesting results which
 improve significantly earlier results obtained by 
 several authors (cf. [Cu-Dr-Lu],[Ko 2], [Ze 1], [Ze 2]).

 Indeed, first we easily deduce new estimates on the 
 relative volume with respect to balls in a generic 
 subspace of $\C^n$ of the 
 plurisubharmonic lemniscates associated to 
 the Lelong class of \psh functions with 
 logarithmic singularities at infinity on $\C^n$ 
 as well as the Cegrell class of \psh functions with 
 bounded Monge-Amp\`ere mass on a bounded
 hyperconvex domain of  $\C^n$.
 
 Then we give estimates on global uniform 
 integrability of the Lelong class of \psh functions 
 with logarithmic singularities at infinity on $\C^n$
 with respect to the Lebesgue measure on any 
 generic subspace. These estimates
 can be considered as precise quantitative versions 
 for the Lelong class of the well known 
 John-Nirenberg inequalities 
 for $\mrm{BMO}-$functions on $\R^n$ (cf.[St]).
 
 In particular we prove that restrictions to any 
 generic subspace $\mb G \sub \C^n$ of plurisubharmonic 
 functions with logarithmic singularities at infinity on 
$\C^n$  are in $\mrm{BMO} (\mb G)$ with a uniform 
  explicit bound on
 their $\mrm{BMO (\mb G)}-$norms depending 
 only on the dimension of 
 $\mb G.$

 Finally we give a general sufficient condition 
 for uniform integrability of a given class of 
 \psh functions on some domain in terms of the
 behaviour of the relative Monge-Amp\`ere 
 capacity of their sublevel sets with respect to this 
 domain. In particular, we deduce a new global 
 uniform integrability result 
 for the Cegrell class of \psh functions of uniformly
 bounded Monge-Amp\`ere masses on a bounded
 hyperconvex domain. 
 \section{Preliminaries}
 Let us recall the classical Polya's inequality
(cf. [Ra], [Ts]). For any compact subset
$K  \sub \C,$ 
\begin{equation}
\lambda_2 (K) \leq    \pi \cdot c (K)^2 \ ,
\label{eq:polyaC1}
\end{equation}
 with equality for a disc,
 where $\la_{2}$ is the area measure on 
 $\C = \R^{2}$ and $c (K)$ is the logarithmic 
 capacity of $K.$

 Besides this inequality, there is a corresponding 
 inequality for sets of the real line $\R \sub \C.$
 Namely, for any compact  subset $K \sub \R,$
 \begin{equation}
 \lambda_1 (K) \leq 4 \cdot  c (K)\ , 
 \label{eq:polyaR1}
 \end{equation}
 with equality for an interval,
 where  $\lambda_1 $ is the lenght measure on 
 $\R$.
 
 Recall that the logarithmic capacity $c (K)$  of 
 a compact subset $K \sub \C$ coincides with its 
 Chebychev constant (cf. [Ra], [Ts]), so that the following 
 formula holds 
 $$ c (K) = \inf _{d \geq1} 
 (\inf\{||P||^{1 \slash d}_K \ ;  P \in  
 \dot{\mP}_{d} \}) \ , $$
 where $\dot{\mP}_{d}$ is the set of monic 
 polynomials of degree $d$ and  $ ||P||_K :=\sup_{z \in K} |P(z)|.$

 In  $\C^n$, it is more convenient to consider the 
 following  Chebyshev constant associated to a compact subset $K \sub \C^n$ (cf. [Al-Ta], [Si 2]) 
 \begin{equation}
 {\mrm{T}}_{B} (K) := \inf _{d\ge 1} (\inf 
 \{||P||^{1 \slash d}_K \ ; \ P \in \C [z], 
 \mrm{deg} (P) = d, ||P||_{B} = 1\}), 
        \label{eq:T}
  \end{equation}
 where $B $ is any regular compact subset of 
 $\C^n$ and $||P||_B := \max_{z \in B} 
 |P(z)|.$ 
 
 If $n=1$, it is easy to prove that the two 
 constants $c$ and $T_B$ are equivalent as we 
 shall see below.
 
 The constant defined by (\ref {eq:T}) is related 
 to the pluricomplex Green function with logarithmic 
 singularities at infinity on $\C^n$, which we 
 will recall below. Its definition 
 is based on the usual Lelong
 class of \psh functions of logarithmic 
 growth at infinity on $\C^n$ defined as follows
 \begin{equation}
        \mathcal{L} (\C^n) := \{u \in PSH (\C^n)\ ;
 \ \sup \{u(z) - \log^{+} \vert z \vert; z \in \C^n \} < + 
 \infty\}.
        \label{eq:L}
 \end{equation}
 The global extremal function with 
 logarithmic growth at infinity  associated to a 
 Borel subset $K \Sub \C^n$ is defined by 
 \begin{equation}
        V_K (z) :=\sup \{u (z) ; u  \in \mL (\C^n), \ \
 u{\vert K}\leq 0\}, z \in \C^n
        \label{eq:V}
 \end{equation}
 and its upper semi-continuous regularization
 $V^*_K$ in $\C^n$ is the pluricomplex Green function 
 with logarithmic singularities at infinity 
 associated to $K$ (see [Za], [Si 1]).
 
 It is well known that $V_{K}$ is locally bounded
 on $\C^n$ if and only if $K$ is non pluripolar in 
 $\C^n$ (see [Si 1], [Si 2]).
 
 By a theorem of Siciak ([Si 2]), we know 
 that if $K \sub \C^n$ si a compact set, then
\begin{equation}
\mrm{T}_{B}(K)=\exp (- \max_{B} V^*_K)
\label{eq:S}
\end{equation}
  The formula (\ref{eq:S}) allows us to extend the 
 definition of the set function $T_B (.)$ to Borel 
 subsets of $\C^n.$ Moreover the extended set
 function is a generalized Choquet capacity on any
 bounded domain in $\C^n,$ which is inner regular and
 outer regular (see [Si 2]). 
 The constant $T_B (K)$ will be called here the 
 {\it relative logarithmic capacity} of $K$ 
 with respect to $B$ in $\C^n$. 
 
 It is also well know that the null sets for this capacity 
 are precisely the pluripolar subsets of  $\C^n$ 
 (see [Si 2]).

Thus if $K \sub \C^n$ is non pluripolar then $- \log \mrm{T}_{B}(K)  = 
 \max_{B} V_{K}^{*} (< + \infty)$ is the best constant for 
 which the following
Bernstein-Walsh inequality holds
\begin{equation}
\sup _{B} u \leq \sup _{K}u   - \log 
\mrm{T}_{B}(K), \quad \forall u\in \mL (\C^n).
\label{eq:BW}
\end{equation}
 
 There is another relative capacity defined using 
 the Monge-Amp\`ere operator (see [Be-Ta 1]).
 Here we choose a normalisation  
 of the usual differential operators on $\C^n$ so 
 that
 $$dd^c := \frac{i}{\pi} \partial
 \overline{\partial}.$$
 Let $\W \Sub \C^n$ be an open set and $K \sub 
 \W$ a compact subset. Then the relative  
 Monge-Amp\`ere capacity of the condenser 
 $(K,\W)$ is defined by the formula (see [Be-Ta 1])
\begin{equation}
\mrm{cap} (K ; \W) := \sup \{ \int_{K} 
(dd^c u)^n ; u \in PSH (\W), - 1 \leq u \leq 0 \}.
        \label{eq:Cap}
\end{equation}
This capacity is related to the so called 
{\it plurisubharmonic measure}
associated to the condenser $(K,\W)$ defined by
\begin{equation}
        h_{K} (z) := \sup \{ u (z) ; u \in PSH (\W), u \leq 0,
       u \vert K \leq - 1 \}, z \in \W.
        \label{}
\end{equation}
Then if $\W \Sub \C^n$ is a hyperconvex open 
set and $K \sub \W$ is a compact subset,
it follows from ([Be-Ta 1]) that
\begin{equation}
        \mrm{cap} (K ; \W) = \int_{K} 
        (dd^c h_{K}^*)^n = \int_{\W} (dd^c h_{K}^*)^n .
        \label{}
\end{equation}
We will need the following Alexander and Taylor's
 comparison 
inequality (see [Al-Ta]).
For a fixed bounded domain $\W \Sub \C^n$ and 
a fixed euclidean ball $\B \sub
\C^n$ such that $\W \sub \B,$ 
\begin{equation}
T_{\B} (E) \leq \exp (- \mrm{cap} (E ; \W)^{- 1 \slash n}). 
\label  {eq:AT}
\end{equation}
 for any Borel subset $E \sub \W.$

 We will also need to define the Cegrell class of \psh
 functions. 
 Let $\W \Sub \C^n$ be a hyperconvex open set.
 Denote by $\m F (\W)$ the class of negative \psh 
 functions $\vphi$ on $\W$ such that there exists 
 a decreasing sequence $(\vphi_j)$ of bounded 
 \psh functions on $\W$ with boundary values $0$ which 
 converges to $\vphi$ on $\W$ and satisfies
 $\sup_j \int_\W (dd^c \vphi_j)^n < + \infty.$ 

 By Cegrell ([Ce 2]), for $\vphi \in \m F (\W),$ the
 Monge-Amp\`ere measure
 $(dd^c \vphi)^n$ is a well defined  Borel measure 
 of finite mass 
 on $\W$ as the weak limit of the sequence
 of measures $(dd^c \vphi_j)^n$, where 
 $(\vphi_j)$ is any decreasing sequence converging 
 to $\vphi$ on $\W$ and satisfying all
 the requierements of the definition.
 
 \section{Relative Polya's inequalities}
 Here we want to compare the relative Lebesgue 
 measure on a generic subspace  $\G \sub \C^n$ 
  with respect to a real euclidean 
 ball  in $\mb G$  with the relative logarithmic 
 capacity in $\C^n$ with respect to the same ball.  
 
 First recall some definitions. A real subspace  
$\mb G \sub \C^n$ is said to be a 
{\it generic subspace} of $\C^n$ if
$\mb G  + J \mb G  = \C^n,$ where $J$ is the complex structure on $\C^n.$ We denote by 
$\mb G ^c := 
\mb G \cap J \mb G $ the maximal complex subspace of $\C^n$
contained in $\mb G$ and set $m := 
\mrm{dim}_{\C} \mb G^c,$ which will be called 
the {\it complex
dimension} of $\mb G.$ Then it is clear that 
$ \mrm{dim}_{\R} \mb G = n + m.$  

If $m = 0$ which means that $\mb G^c = (0),$  
the subspace 
$\mb G$ is said to be {\it totally real}. If $m = n$ then
$\mb G = \C^n.$

 It is easy to see that  $\mb G \sub \C^n$ is a 
 generic subspace of complex dimension $m$ if and only if
 there is a unitary automorphism
$U: \C^n \lgra \C^n$ such that $ U (\mb G) = 
\C^m \times 
\R^{n - m} \sub  \C^m \times \C^{n - m} =  \C^n$. 
 
Observe that the subspace  $\mb G \sub \C^n$ 
 is non pluripolar in 
 $\C^n$ precisely when $\mb G$ is a generic subspace.

The subspace $\mb G \sub \C^n$ will be endowed 
with the induced euclidean structure and the 
corresponding Lebesgue measure which will be 
denoted by $\la_{n + m}.$ 

Now we can state our version of Polya's inequality
which is the main result of this section.
\begin{thm} 1) For any complex euclidean closed ball
$\B \sub \C^n$ and any Borel subset $K \sub \B,$  
\begin{equation}
 \frac{\la_{2 n} (K)}{\la_{2 n}(\B)} \ 
\leq \ c_ n \ T_{\B} (K)^{2}.
\label{eq:PolyaC}
\end{equation}
where 
\begin{equation}
 c_n := \frac{4^{n} (n!)^2}{(2n - 1)!}.
\label{eq:cn}
\end{equation}
2) Let $\mb G \sub \C^n$ be a generic 
real subspace of  complex dimension $0 \leq m \leq n - 1.$ 
Then for any real euclidean closed ball $B \sub \mb G$ and 
any Borel subset $K \sub B,$  
\begin{equation}
 \frac{\la_{n + m} (K)}{\la_{n + m} (B)} \ 
\leq \ 8  (n + m) \ T_{B} (K).
\label{eq:PolyaR}
\end{equation}
\end{thm}
For the proof of relative Polya's inequalities, we start to look to the simplest case
where $n = 1.$ 
\begin {lem} 1) For any closed disc $\mb D \sub \C$
and any Borel subset $K \sub \mb D,$
\begin{equation}
 \frac{\la_{2} (K)}{\la_{2} (\mb D)} \ \leq \ 4 \ T_{\mb D}  (K)^{2}.
 \label{eq:PolyaC1}
\end{equation}
2) For any real closed intervall $\mb I \sub \R$
and any Borel subset $K \sub \mb I$
\begin{equation}
 \frac{\la_{1} (K)}{\la_{1} (\mb I)} \ \leq \ 4 \ T_{\mb I} 
(K).
  \label{eq:PolyaR1}
\end{equation}
\end{lem}
We don't know if $4$ is the best constant in these inequalities.\\
\demo 
1) By regularity of the Lebesgue measure and the relative 
logarithmic capacity in $\C$, we can assume that $K$ is  a
non polar compact subset.  We can also assume that
$\C \sm K$ is connected since $\la_{2} (K) \leq
\la_{2} (\hat{K})$ and $T_{\mb D} (K) = T_{\mb D} (\hat{K}).$
 Then the extremal function $V_{K}^{*}$ is a subharmonic
 function on $\C$ which coincides with the Green function of
 $\C \sm K$ with a pole at infinity. Therefore it can be represented by the formula
 $$ V_K^{*} (z) = \int_{K} \log |z - \zeta| d \mu (\zeta) - 
 \log c(K), \qquad \forall z \in \C,
 $$
 where $\mu := (1 \slash 2\pi) \Delta V_{K}^{*}$ is 
 the normalized equilibrium 
 measure of $K.$ 
 From this representation formula,  we get the 
 estimate
 $$ \max_{\mb D} V_{K}^{*} \leq \log (2 R) - 
 \log c (K)\ , $$ 
 where  $R$ is the radius of the disc  $\mb D \sub 
 \C$. This inequality implies that  
\begin{equation}
 c (K) \ \leq \  2 R \ \mathrm{T}_{\mb D}(K). 
 \label{eq:comp-cT} 
 \end{equation}
 Therefore using the inequality (\ref{eq:polyaC1}),
 we get from (\ref{eq:comp-cT}) the estimate
 $$\lambda_2 (K) \ \leq \ 4  \ \lambda_2 (\mb D) \
\mathrm{T}_{\mb D} (K)^{2},$$ 
 which is the required estimate. 

 2) In the real case we prove in the same way
 that
  $$c (K) \ \leq  \ 2 R \ \mathrm{T}_{\mb I}(K),$$
 where $R$ is the radius of the interval $\mb I.$
 Therefore using the inequality (\ref{eq:polyaR1}), 
 we get 
 $$\lambda_1 (K) \ \leq  \ 4  \ \lambda_1 (\mb I) \ 
 \mathrm{T}_{\mb I} (K),$$
 which is the required inequality.
 \fin 
 To prove our theorem, we need the following elementary 
slicing lemma.
\begin{lem}
1) Let  $\B \sub \C^n $ be any complex euclidean closed ball, $K \sub \B$
a Lebesgue measurable subset and  $ a \in
\partial \B$. Then there exists a complex line 
$L_a \sub \C^n$ passing through the point 
$a$ such that $\la_2 (\B \cap L_a) > 0$
 and
 \begin{equation}
    \frac{\la_{2 n} (K)}{\la_{2 n} (\B)} \ \leq \
    c^{\prime}_n  \ \frac{\la_2 (K \cap
 L_a)}{\la_2 (\B \cap L_a)}.
 \label{eq:sliceC}
 \end{equation}
where $ c^{\prime}_n = c_n/4 =\frac{ 4^{n-1} (n!)^2}{(2n - 1)!}.$\\ 
2) Let  $B \sub \R^{N}$ any
 euclidean ball, $K \sub B$ any Lebesgue 
 measurable subset and  $a \in B.$
 Then there exists a real line $l_a \sub \R^{N}$
 passing through the point $a$ such that 
 $\la_1 (B \cap l_a) > 0$ and
\begin{equation}
    \frac{\la_{N} (K)}{\la_{N} (B)} \ \leq \  
   2 N  \ \frac{\la_1 (K \cap
l_a)}{\la_1 (B \cap l_a)}.
\label{eq:sliceR}
\end{equation}
\end{lem}
 Observe that $c_n \sim 2\sqrt{\pi}n^{3/2}$ as $n \to \infty.$, we conjecture
 that the inequality (\ref{eq:sliceC}) is true with the constant  $
 c^{\prime}_n = n$. The inequality (\ref{eq:sliceR}) could be  deduced from
 ([BG], lemma 3) with the constant $N$ but
 the proof given there is not clear for us. So we decided to give 
 proof which uses the same idea of symetrization but leads to
 the constant $2 N$ instead of $N$, unless the point $a$ in the lemma coincides with the
 center of the ball $B.$ \\
\demo 1) We can of course assume that $n \geq 2.$
Since our inequality is invariant under 
translation, we can also assume that
$a = 0 \in \partial \B$ is the origin and $\la_{2 n} (K) > 0.$ \\
Now assume by contradiction that the inequality 
(\ref{eq:sliceC}) is not true. Then we will have 
\begin{equation}
 \lambda_2 (K\cap L ) < \frac{\lambda _{2 n} (K)} 
{c^{\prime}_n \ \lambda_{2 n} (\B)} \ \lambda_2 (\B \cap L),
\label{eq:sliceC'}
 \end{equation} 
 for any complex line $L$ passing through the 
 origin $a=0$ such that $ \lambda_2 (\B \cap L) > 0.$\\
Since relative volume and relative area are invariant under
non singular affine transformations,
 we can assume that $\B = \{z = (z_1, z_2,\cdots,
 z_n) \in \C^{n} \ ; \ \ \vert z_1 -
R\vert^2 + \vert z_2\vert^2 + \cdots +
\vert z_n\vert^2 < R^2\}$ and
$L_w =\{ \zeta .w \ ;\ \zeta \in \C \} $
where $ w =(w_1,\ldots,w_n)\in S^{2 n - 1}.$
Then  $L_w \bigcap B =\{ \zeta \cdot w \ ;\  \vert \zeta
\vert^{2} < 2 R \ \Re (\zeta w_1) \}$ is  the  disc centred at $R w_1$ of radius
$   R \ \vert w_1\vert$ 
which by
the last inequality  leads to
\begin{equation}
\lambda_2 (K\cap L_w ) < \frac{\lambda _{2 n} (K)}
{c^{\prime}_n \ \lambda_{2 n}
(B)} \  \pi  R^{2} \ \vert w_1\vert^{2},\forall w \in
 S^{2 n - 1},
\label{eq:area}
 \end{equation}
 Now, integrating in polar coordinates and using the invariance of the sphere
$S^{2 n - 1}$ by rotation,  we
obtain the formula
\begin{eqnarray*}
\lambda _{2 n} (K) &=&   \frac{1}{2 \pi} \int _{S^{2 n - 1}} \int _{\vert \zeta
\vert < 2 \ R \ \vert \zeta w_1\vert }
\vert \zeta \vert^{2 n - 2}
\chi _{K}(\zeta \cdot w) d \lambda_{2} (\zeta)
d \sigma _{2 n - 1} (w)
  \\
  &\leq &
 \frac{ {2^{2 n - 2} R^{2 n - 2}}}{2\pi}  \int_{S^{2 n - 1}}
  \vert w_1\vert^{2 n - 2}
\int _{\vert \zeta \vert \leq 2 R  \vert  w_{1}\vert}
\chi_{K} (\zeta \cdot w) d \lambda_{2} (\zeta)
 d \sigma _{2 n - 1}(w),
\end{eqnarray*}
where $\chi_{K}$ is the characteristic function of
the set $K.$ \\
Using inequality (\ref{eq:area}), we deduce from
the last inequality that
 \begin{equation}
 \lambda_{2 n} (K) <
   {2^{2 n-2} R^{2 n}} \frac{\la_{2 n}(K)} {2 c^{\prime}_n \  \
\lambda_{2 n}(B)}
 \int_{S^{2 n - 1}} \vert w_1\vert^{2 n} d
 \sigma _{2 n - 1}(w)
 \label{eq:polar2}
 \end{equation}
 Now, an elementary computation using spherical
 coordinates leads to the formula 
 \begin{equation}
 \int_{S^{2n - 1}} \vert w_1 \vert ^{2n} d
 \sigma _{2n - 1}(w)
 =\frac{4n (n!)^2}{(2n)!}\tau_{2n}
\label{eq:for}
 \end{equation}
 where  $\tau_{2n}$ is the volume of the euclidean
 unit ball in $\R^{2n}.$\\
 The last formula (\ref{eq:for}) combined with
 (\ref{eq:polar2}) leads finally to the inequality
 $$ \lambda_{2 n} (K) <
   \frac{ {2^{2 n-2}R^{2n}}\la_{2 n}(K)} {2c^{\prime}_n \
\lambda_{2 n}(B)} 4n \frac{(n!)^2}{(2n!)}
 \tau_{2 n}  =\lambda_{2 n} (K),$$
which yields a contradiction. \\
 2)  As in the complex case, we assume that
$a = 0$ is the origin, $\lambda_{N} (K) > 0$  and the ball $B$ of
 radius $1$. \\
First, observe that $\lambda_1(B \cap l_a)\le 2 $ for any real
  line $l_a$ passing through the point $a$ , then to show (\ref{eq:sliceR})
 it is enough to prove that 
$$
\frac{1}{N}\frac{\la_{N} (K)} {\la_{N} (B)} \le \lambda_1(K \cap l_a) $$
for some real line $l_a.$\\
Assume by contradiction that the last inequality
 is not true. Then we will have
\begin{equation}
 \lambda_1 (K\cap l ) < \frac{1}{N} \frac{\la_{N} (K)} {\la_{N} (B)} 
\label{absR}
 \end{equation}
 for any real line $l$ passing through the
 origin $a=0.$
 \\
 Let $\tilde {K}$ be the annulus with the same  center $x_0$  as  $B$ and of
  radii
 $r$ and $ 1$ \ $ (r<1)$
such that $ \la_{N}(\tilde{K})=\la_{N}(K)$\\
then $$
r = \Bigl(1 - \frac{\la_{N}(K)}{\la_{N}(B)}\Bigr)^
{1 \slash N}.
$$
Denote by $e(\tilde{K}) := 1 - r$ the depth of the annulus
 $\tilde{K}$,
 then
\begin{eqnarray*}
e(\tilde{K}) &= &  1 - r \\
&=&  1- \Bigl(1-\frac{\la_{N}(K)}{\la_{N}(B)}\Bigr)^{\frac{1}{N}} \\
&\ge &  \frac{1}{N} \frac{\la_{N}(K)}{\la_{N}(B)}.
 \end{eqnarray*}
The last inequality together with (\ref{absR}) lead to
\begin{equation}
e(\tilde{K}) >  \la_1( {K}\cap l)\quad 
\label{abR}
\end{equation}
\mbox{for any  real
 line $l$ passing through}\ \ $ a.$\\
 Now, observe that, if $l$ any real line passing through the origin such that 
$l\cap B(x_0, r) \ne \es$,    then 
$\la_1((\tilde {K}\cap l) \ge 2  e(\tilde{K})$ 
   and hence from (\ref{abR})
   we derive the inequality 
\begin{equation}
\la_1(\tilde {K}\cap l) > 2  \la_1( {K}\cap l),
\label{incR}
\end{equation}
for  all real line passing through the origin $a=0$. \\
Now, following ([Br-Ga]), we  construct a  set $K^{(s)}$ in the following way:
On each real line $l$ passing through the point $a=0$, we choose the best far
 segment of 
$\tilde{K}\cap l$ of  length  $\la_1(K\cap l).$\\
Then from the inequality (\ref{incR},) we get $K^{(s)}\subset \tilde{K}$
 and therefore
\begin{equation}
\la_{N}(K^{(s)})< \la_{N}(\tilde{K}).
\label{sumR}
\end{equation}
On the other hand, by the construction of the set  $K^{(s)}$,  if $ \tau \in
 K \cap
l \backslash  K^{(s)} $ and $ t \in ( {K^{(s)}} \cap l) \backslash K $
then $ |\tau |\le |t|$  and 
 since 
$\la_1 ({K\cap l  })=\la_1( {K^{(s)}} \cap l )
 $
 then
$$\int_{K\cap l  } |\tau |^{N-1} d\tau \le
\int _{ K^{(s)}\cap l} |t|^{N-1} dt.$$
 Now, integrating in polar coordinates and using  the last inequality,  we
obtain
\begin{eqnarray*}
\lambda _{N} (K) &=& \frac{1}{2}   \int _{S^{N - 1}} \int _{\R }
\vert \tau  \vert^{N-1}
\chi _{K}(\tau \cdot w) d\tau 
d \sigma _{N- 1} (w)
  \\
 &=&  \frac{1}{2}   \int _{S^{N - 1}} \int _{K\cap l_w }
\vert \tau  \vert^{N -1}
 d\tau 
d \sigma _{N - 1} (w)\\
 &\le&   \frac{1}{2}   \int _{S^{N - 1}} \int _{K^{(s)}\cap l_w }
\vert t \vert^{N -1}
 dt
d \sigma _{N - 1} (w)\\
  &\le &  \frac{1}{2} 
  \int _{S^{N - 1}} \int _{\R }
\vert t \vert^{N -1}
\chi _{K^{(s)}}(t \cdot w) dt
d \sigma _{N - 1} (w)\\
&\le & \lambda _{N} (K^{(s)})
\end{eqnarray*}
where $\chi_{K}$ is the characteristic function of
the set $K$ and $l_w=\{t\cdot w\ ; \ t\in \R \},$ \\
which contradicts the inequality (\ref{sumR}).\fin
 Now we are ready for the proof of the Theorem. \\
 \underline {Proof of the theorem:} 
 1) By interior and exterior
 regularity of the Lebesgue measure and the relative logarithmic
 capacity, we can assume that $K \sub \B$ is a
 compact set of non empty interior in $\C^n$
 so that $\la_{2 n} (K) > 0$ and $T_{\B} (K) > 0.$
Moreover, considering $\vep-$neighbourhoods of $K$ in $\C^n$, we can assume
that $K$ is regular in the sense that $V_K$ is continuous on $\C^n.$ Therefore  
 $V_K \in \mathcal{L} (\C^n)$ and there exists  $a\in \partial \B$
 such that $V_K (a) = \sup_{\B} V_K.$ By translation, we can assume  that $a = 0$ is the origin in $\C^n.$ \\
 By the complex slicing lemma, there exists a complex line
 $L \sub \C^n$ passing through the point $a$ such 
 that $\lambda_2 (K \cap L) > 0$ and
 \begin{equation}
\frac {\lambda_{2n} (K)}{\lambda_{2n} (\B)} \leq 
  {c^{\prime}_n}  \ \frac{\lambda_2 (K \cap L)}{\lambda_2 (\B \cap L)}.
\label{eq:rel-ineq1}
 \end{equation}
  Since  $a \in L$ and $V_K (a)= \max_{\B} V_K$, it follows
 that $T_{\B \cap L} (K \cap L) \leq T_\B (K)$ and then from
  (\ref{eq:rel-ineq1})and (\ref{eq:PolyaC1}) we deduce that 
 \begin{equation}
\frac {\lambda_{2n} (K)}{\lambda_{2n} (\B)} \leq 4 c^{\prime}_n T_\B (K)^2,
 \label{eq:BWL1}
\end{equation}
which is exactly the required inequality (\ref{eq:PolyaC}). \\
2) Now assume that $\mb G \neq \C^n$ is a generic subspace of
 complex dimension $1 \leq m \leq n - 1$ 
 (the totally real  case $m = 0$ is treated in the same way).
 By the invariance of the Lebesgue measure and the relative
 capacity $T_ B$ by unitary transformations, 
 we can assume that  $\mb G = \C^m \times \R^{n - m}.$
 By regularity properties of the Lebesgue measure 
and the relative capacity $T_B$, we can assume that
 $K \sub B$ is a compact subset of non 
 empty interior in $\mb G$ so that 
 $\la_{n + m} (K) > 0.$ Let us prove that $T_B (K) > 0.$
 Indeed, since $K$ is a compact subset of non 
 empty interior in $\mb G$,
 there exists an intervall $I \sub \R$ of positive lenght and a
 disc $\mb D \sub \C$ of positive radius such that 
 $ \mb D^m \times I^{n - m} \sub K.$ 
 Then by the
 product property of the extremal function 
 (cf. [Si 1]), we get 
 $$V_K (z,\zeta) \leq
  \max_{1 \leq i \leq p, 1 \leq  j \leq n - m} \{V_{\mb D} 
  (z_{i}), V_I (\zeta_j) \}, $$
   for any $z = (z_{1}, \ldots , z_{m}) \in \C^m$ and
   $\zeta \in \C^{n - m}.$ 
Therefore $V_{K}$ is locally bounded on $\C^n$
 and then $T_B (K) > 0.$ Considering $\vep-$neighbourhoods of $K$ in $\G,$ we
 can assume by regularity that $K$ is a regular compact set in the sens that $V_K$ is continuous in $\C^n.$ \\
 Then   $V_K \in \mathcal{L} (\C^n)$ and there exists  $a \in 
 B$ such that $V_K (a)= \sup_{B} V_K.$ \\
 By translation we may assume that $a = 0$ is the origin
 in $\mb G.$
 Then by the real slicing lemma, there exists a real line 
 $l \sub \mb G$ passing through the point $a = 0$ 
 such that $\lambda_1 (K \cap l) > 0$ and
 \begin{equation}
 \frac{\lambda_{n + m} (K)}{\lambda_{n + m} (B)} \leq  
 \ 2(n + m) \ \frac{\lambda_1 (K \cap l)}{\lambda_1 (B \cap l)}.
\label{eq:rel-ineq}
 \end{equation}
 Let $L := l + i \cdot l$ be the complex line in 
 $\C^n$ generated by the real line $l.$ Since  
 $a = 0 \in l$ and $V_K (a)= \sup_{B} V_K$,  
 it follows that $T_{B \cap l} (K \cap l) \leq T_B (K)$ and
 then from (\ref {eq:PolyaR1}) and (\ref{eq:rel-ineq}) we
 deduce that  \begin{equation}
 \frac{\lambda_{n + m} (K)}{\lambda_{n + m} (B)}
 \leq 8 (n + m) T_B (K),
 \label{eq:BWL}
\end{equation}
 which is exactly the required inequality (\ref{eq:PolyaR}).
\fin  
 It is interesting to observe that from the formula (\ref{eq:S}) it follows that  our relative Polya's
 inequalities leads to the following quantitative
 versions of Bernstein-Walsh inequalities.
 \begin{cor}
 1)  For any closed complex euclidean ball $\B \sub 
 \C^n$, any Borel subset $K \sub \B$ and any function  $u \in \mathcal{L}( \C^n),$
 \begin{equation}
\sup _{\B} u \ \leq \ \sup _{K} u + \frac{1}{2}
\log ({c_n}) - \frac{1}{2} \log  
\frac{\la_{2n} (K)}{\la_{2 n} (\B)},
\end{equation}
where $c_n$ is the constant given by (\ref{eq:cn}).\\
2) Let $\mb G \sub \C^n$ be any generic subspace
of complex dimension $m \leq n - 1.$ Then for any closed
real euclidean ball $B \sub \mb G,$ 
any Borel subset $K \sub B$ and  any function  $u \in \mathcal{L}(\C^n),$
\begin{equation}
\sup _{B} u  \ \leq  \ \sup _{K} u + \log \bigl( 8 (n + m)\bigr) - \log \frac{\lambda_{n + m} (K)}{\lambda_{n + m} (B)}.
 \label{eq:BWreal}
\end{equation}
\end{cor}
 Let us mention that in the totally real case 
$\mb G = \R^n,$ inequalities like (\ref{eq:BWreal}) 
where obtained earlier by A. Brudnyi (cf. [B 1]).

 From relative Polya's inequalities (\ref{eq:PolyaC}),
 (\ref{eq:PolyaR}) and Alexander-Taylor's
 inequality (\ref{eq:AT}), we deduce the following interesting
 comparison inequalities between relative volumes and the relative
 Monge-Amp\`ere capacity. These inequalities show that the Lebesgue measure on any generic
subspace of a hyperconvex domain $\W \Sub \C^n$ is dominated by
capacity in a strong sense and then by a result of S. Kolodziej, it belongs the image of the
complex Monge-Amp\`ere operator acting on the class of bounded \psh functions
 on $\W$ (see [Ko 1], [Ko 2], [Ce 1]).
\begin{cor} 1)For any complex euclidean ball $\B \sub \C^n$
 and any Borel subset $ K \sub \B,$
 \begin{equation}
\frac{\lambda_{2 n} (K)}{\lambda_{2 n}(\B)}
\leq c_n \exp (- 2 \ {\mrm{cap}(K;\B)}^{- 1 \slash n}),
\label{eq:MAC}
\end{equation}
where $c_n$ is the constant given by (\ref{eq:cn}). \\
2) Let $\mb G \sub \C^n$ be a generic 
real subspace of  complex dimension $0 \leq m \leq n - 1.$ 
Then for any euclidean ball $B \sub \G$ and 
any Borel subset $K \sub B,$  
\begin{equation}
 \frac{\la_{n + m} (K)}{\la_{n + m} (B)} \ 
\leq \ 8  (1 + \sqrt{2}) \ (n + m) \
 \exp (- \ {\mrm{cap}(K;\B)}^{- 1
\slash n}),
\label{eq:MAR}
\end{equation}
 where $\B$ is the euclidean ball in $\C^n$ such that $\B \cap
 \G = B$.
\end{cor}

 \demo 1) The inequality (\ref{eq:MAC}) is a direct 
 consequence of (\ref{eq:AT}) and (\ref{eq:PolyaC}).\\
 2) Let us prove the inequality (\ref{eq:MAR}). Since both the
 relative volume and the relative capacity are invariant
under non singular  affine transformations, we can assume that
 $\mb G = \C^m \times \R^{n - m}, B$ is the unit real euclidean
 ball in $\mb G$ and $\B$ is the unit complex euclidean ball in $\C^n.$
 Then by  (\ref{eq:PolyaR}), we have
\begin{equation} \frac{\la_{n + m} (K)}{\la_{n + m} (B)} \ 
\leq 8 (n + m) T_B (K).
\label{eq:star}
\end{equation}
On the other hand, by (\ref{eq:AT}), we have
$$ T_{\B} (K) \leq  \exp (- {\mrm{cap} (K;\B)}^{- 1
\slash n}).
$$
 So to prove the inequality (\ref{eq:MAR}), it remains to
 estimate $T_{B} (K)$ from above by $T_{\B} (K).$
 Indeed, from the definition of the extremal function $V_B,$
 it follows that
$$ V_K (z) \leq \max_B V_K + V_B (z), \forall z \in \C^n.$$
Therefore, we get 
\begin{equation}
 T_B (K) \leq e^{\max_\B V_B} \ T_\B (K).
\label{eq:TB}
\end{equation}
It remains to estimate $\max_\B V_B.$
 Since $\R^n \sub \mb G,$ the euclidean unit ball $B$ in
 $\mb G,$  contains the  euclidean unit ball
 $D$ of $\R^n$ and then $V_B \leq V_D$ on $\C^n,$ which implies  that $ \max_\B
 V_B \leq \max_\B V_D.$
 Now by Lundin's formula (cf. [Lu],[Sa 2], [Kl]), we have
\begin{equation}
 V_D (z) = \max \{ \log \vert h (\xi \cdot z)\vert \ ; \
 \xi \in S^{n - 1} \},  z \in \C^n,
\label {eq:Lu}
\end{equation}
 where $ h (\zeta) := \zeta + \sqrt{\zeta^2 - 1}$
 for $\zeta \in \C$, with the right branch of the square root,
$S^{n - 1} =  \partial D$ is the euclidean unit sphere of $\R^n \sub \C^n$ and 
 $\xi \cdot z = \sum_{1 \leq j \leq n} \xi_j \cdot z_j$.
 It is easy to see from the formula (\ref{eq:Lu}) that
 $$\max_\B V_D = \max_{\vert z \vert = 1} V_D (z)= \max_{\vert \zeta \vert
 = 1} \log \vert h (\zeta)\vert =  \log (1 + \sqrt{2})$$
 and then $\exp (\max_\B V_B) \leq \exp (\max_\B V_D) =  1 + \sqrt{2},$
 which by the inequality (\ref{eq:TB})and (\ref{eq:star})
 implies the
 required inequality (\ref{eq:MAR}).
\fin 
\noindent {\bf Remarks:} 1) Polya's inequalities 
(\ref{eq:PolyaC}) and
(\ref{eq:PolyaR}) can be stated
in one formula as follows.
Given a generic 
subspace $\G \sub \C^n$ 
of complex dimension $0 \leq m \leq n,$ then 
for any euclidean ball
$B \sub \G $ and any Borel subset $K \sub B,$ we have
\begin{equation} \frac{\la_{n + m} (K)}
{\la_{n + m} (B)} 
\leq c_{n,m} 
T_{B} (K)^{1 + [m \slash n]},
\label{eq:PG}
\end{equation}
where $c_{n,m} := 8 (n + m)$ if $0 \leq m \leq n - 1$
and $c_{n,n} := c_n.$

We can deduce from the general relative Polya's 
inequality (\ref{eq:PG})  analoguous inequalities 
in terms of relative volume and relative logarthmic capacity with respect to 
balls associated to any fixed real norm on the generic space
$\mb G.$ Indeed, if we denote by $\vert . \vert$ 
the euclidean norm and we are given another real norm 
$\Vert . \Vert$  on $\mb G,$ then there exists two 
constants $\alpha, \beta > 0$ such that
$$ \alpha \Vert z \Vert \leq 
 \vert  z \vert \leq  \beta \Vert z\Vert, \forall 
z \in \mb G.$$
Then given a ball $B'$ for the norm $\Vert . \Vert$,
there exists a ball $B$ for the norm $\vert . \vert$ 
such that
$\alpha \cdot B \sub B' \sub \beta \cdot B.$
Then it follows easily from (\ref{eq:PG}) that
for any Borel set $K \sub B',$ we have
\begin{equation} \frac{\la_{n + m} (K)}{\la_{n + m} (B')} 
\leq c_{n,m} (\beta \slash \alpha)^{n + m} 
T_{B'} (K)^{1 + [m \slash n]}.
\label{eq:PGbis}
\end{equation}
2) Observe that relative Polya's inequalities proved above are optimal 
as far as the exponents are concerned.
Indeed we will use inequality (\ref{eq:PGbis}) for the sup-norm, since in
 this case, explicit computations can be made using  the product formula for the
relative 
logarithmic capacity. Let  $B_1, ..., B_n$ be
regular sets in $\C$, $K_1, ..., K_n$ Borel subsets such that
$K_j \sub B_j$ for $j = 1, ..., n$ and set  $K :=
 K_1 \times ... \times K_n$ and $B := B_1  \times ... \times B_n$.
Then using the product property for the extremal function (cf. [Si 1]), we get
the formula
\begin{equation}
T_B (K) = \min_{1 \leq j \leq n}  \{ T_{B_j} (K_{j})\}.
\end{equation}
In the case where $\mb G = \C^n,$ take $B'$ to be 
 the closed unit polydisc $\Delta^n$ in $\C^n$ and 
 $K_{r} := \{z \in \Delta^n ; \vert z_{1} \vert \leq r\}.$ 
 Then the relative volume of
 $K_r$ with respect to $\Delta^n$ is $\la_{2 n} (K_r) \slash
\la_{2 n}(\Delta^n)=  r^{2}$ while its relative
 logarithmic capacity is $T_{\Delta^n} (K_{r}) = r.$
 Then by (\ref{eq:PGbis}) this prove that the exponent $2$
 in the complex Polya's inequality (\ref{eq:PolyaC}) is the
 best possible.

 In the totally real case, we can assume that
 $\mb G = \R^n$ and consider an analoguous example 
 with intervals. Take $B'$ to be the unit $n-$cube 
 $\mb I^n,$ where $\mb I := [- 1 , + 1]$ is the closed unit real
 interval and
 define $I^n (r) := \{x \in \mb I^n ; \vert x_{1}\vert  \leq r\}.$
 Then it is easy to see that
 $$T_{\mb I^n} (I^n (r)) = \frac{r}{1 + \sqrt{1 - r^2}} \sim \ \frac{r}{2}, \ \
 \mrm{as} \ r \to 0,$$
 while the relative $n-$volume of $I^n (r)$ with respect to
 $\mb I^n$
 is equal to $r$, which proves by (\ref{eq:PGbis}) that the
 exponent $1$ in
 Polya's inequality
 (\ref{eq:PolyaR}) is the best possible in this case.

 Now if $\mb G = \C^m \times \R^{n - m}$ with
 $(1 \leq m \leq n - 1,$ it is enough to take
 $B' = \Delta^m \times \mb I^{n - m}$ and $K_{r} := 
 \Delta^m \times I^{n - m} (r).$
 Then $T_{B'} (K_r) \sim r \slash 2$ as $r \to 0$,
 while $\la_{n + m} (K_r)\slash \la_{n + m} (B') = r,$
 which prove again
 by (\ref{eq:PGbis}) that the exponent $1$ in Polya's
 inequality (\ref{eq:PolyaR}) is the best possible in this
 case. 
\section{Relative size of \psh lemniscates}
 Here we want to deduce from relative Polya's inequalities 
 an estimate on the relative size of plurisubharmonic
 lemniscates associated to two important classes of \psh
 functions.
 Let us start with estimating precisely the size of the 
 lemniscates associated to the Lelong class $\mL (\C^n).$
 \begin{thm}
 1) For any complex euclidean closed ball $\B \sub \C^n$ 
 and any  $u \in \mathcal{L} (\C^n)$ with
$\max_{\B} u = 0,$ 
\begin{equation}
\frac{\lambda _{2n} (\{z \in \B  ;   u(z) \leq - s\})}{\lambda_{2n} (\B)}
\leq {c_n} 
 \ e^{- 2 s}, \ \ \forall s > 0,
 \label{eq:lemC}
\end{equation}
where $c_n$ is the constant given by (\ref{eq:cn}).\\
2) Let $\G \sub \C^n$ be a generic real subspace of complex dimension
$m \leq n - 1.$ Then for any real euclidean closed ball 
$B \sub \G$ and
any  $u \in \mathcal{L}(\C^n)$ with
$\max_{B} u = 0,$
\begin{equation}
\frac{\la_{n + m} (\{x \in B  ;\quad  u(z) \leq - s\})}{ \la_{n + m} (B)} \ \leq \ 8 \ (n + m) \
  e^{- s},  \ \ \forall s > 0.
  \label{eq:lemR}
\end{equation}
  \end{thm}
\demo  1) Let $\B \sub \C^n$ be an arbitrary complex ball
and $u \in \m L (\C^n)$ with $\max_{\B} 
u = 0.$ Set $ E_t (u) := \{z \in \B ; u (z) \leq t\}$ 
for
$t < 0.$ Then 
$u - t \leq V_{E_t (u)}$ on $\C^n$ and then
 $- t = \max_{\B} u - t \leq \max_\B V_{E_t(u)}$. 
 This implies that
$T_\B (E_t(u)) \leq e^t$ for any $t < 0.$
Now in order to get
the estimate (\ref{eq:lemC}), it is enough to apply the complex Polya's inequality 
(\ref {eq:PolyaC}) to the Borel set $E_t(u)$ with
 $s = - t$.  To prove the estimate
(\ref{eq:lemR}), we proceed in the same way using
the real Polya's inequality (\ref{eq:PolyaR}).
\fin
Observe that estimates of \psh lemniscates were 
obtained in the complex case earlier by the 
third author in a more general context
but with less precise exponents (cf. [Ze 2], [Ze 3]).

In particular, observing that $(1 \slash d) \log \vert P\vert \in \m L (\C^n)$
for any polynomial $P \in \C [z]$ with degree $d \geq 1,$ we obtain the
 following precise estimate for polynomial lemniscates.
\begin{cor}
1) For any complex ball $\B \sub \C^n$ and  
any polynomial  $P \in \C [z] $ of degree 
$d \geq 1$ satisfying  $||P||_{\B} = 1,$ we have
\begin{equation}
\frac{\lambda _{2n} (\{z \in \B ;\quad  |P (z)| \leq 
\varepsilon ^d \})}{ \la_{2n} (\B)} \leq {c_n}  \
\varepsilon^{2}, \ \forall \varepsilon \in ]0 , 1],
\end{equation}
where $c_n$ is the constant given by (\ref{eq:cn}).\\
2) Let $\G \sub \C^n$ be a generic subspace of 
complex dimension $0 \leq m \leq n - 1.$ Then
for any real euclidean ball $B \sub \G$ 
and  any polynomial  $P \in \C [z] $ of degree 
$d \geq 1$ satisfying  $||P||_{B} = 1,$  the following
 estimate holds
\begin{equation}
\frac{\la_{n + m} (\{z \in B ;  |P (z)| \leq
 \varepsilon^d \})}{ \la_{n + m} (B)} \ \leq \
 {8 (n + m)} \ \varepsilon, 
\ \forall \vep \in ]0 , 1].
\end{equation}
\end{cor}
All these estimates are optimal as far as the 
exponents are concerned (see Remarks above).
 The first inequality is an improvement of 
 previous results (see [Cu-Dr-Lu], [Ze 2], [Ze 3]) and 
 answers a question
 asked by the third author in ([Ze 2]).  
 In the totally real case where $\G = \R^n,$
 the second inequality appears also in ([Br-Ga]).
 
 Now let us estimate the size of \psh lemniscates 
 associated to the Cegrell class $\m F (\W).$ 
\begin{thm} Let $\W \Sub \C^n$ be a hyperconvex 
open set. Then
for any \psh function $\vphi \in \m F (\W)$ with 
$\int_\W (dd^c \vphi)^n   \leq 1,$ we have
\begin{equation} \la_{2 n} (\{z \in \W ; \vphi (z) \leq - s
 \}) \leq c_n \tau_{2 n} (\W)
e^{- 2 s}, \ \ \forall s > 0,
\label{eq:CegC}
\end{equation}
where $ \tau_{2 n} (\W)$ is the volume of the smallest euclidean ball
of $\C^n$ containing $\W$ and $c_n$ is the constant given by (\ref{eq:cn}).

Moreover, if  $\G \sub \C^n$ is a generic subspace
of complex dimension $m \leq n - 1$ such that $D := 
\W \cap \G \neq \es,$ then for any $s > 0,$
\begin{equation}
\la_{n + m} (\{z \in D ; \vphi (z) \leq - s \}) \ \leq 
\ 8  (1 + \sqrt{2}) \ (n + m)  \
 \tau_{n + m} (D) \ e^{- s},
  \label {eq:CegR}
\end{equation}
 where $ \tau_{n + m} (D)$ is the volume of the smallest
 euclidean ball of $\G$ containing $D.$
\end{thm}
For the proof of this theorem, we will need the 
following elementary lemma.
  \begin{lem} Let $\W \Sub \C^n$ be a hyperconvex open set.
Then for any $\vphi \in \mF (\W),$
\begin{equation} 
  \mrm{cap}(\{z \in \W ; \vphi (z) \leq - s \};\W) \leq s^{- n} 
  \int_{\W} (dd^c \vphi)^n, \ \ \forall s > 0.
\label{eq:F}
\end{equation}
  \end{lem}
  \demo 1) Assume first that $\vphi$ is a bounded \psh
 function on $\W$ with boundary values $0$ and finite
 Monge-Amp\`ere mass on $\W.$
  Let $s > 0$ be fixed and $K \sub \W(\vphi;s):= \{z \in \W ; \vphi (z) \leq - s \}$
  any fixed regular compact set in the sense that
  the plurisubhamonic measure
  $h_{K}$ the condenser $(K,\W)$ is continuous 
  on $\W.$ 
 Since $h_{K}$ and $\vphi$ have boundary 
 values $0,$ from the comparison principle
 (see [Be-Ta 1], [Kl]) it
 follows that
  $$\mrm{cap} (K;\W) = \int_{K} (dd^c h_{K})^n
  \leq \int_{\{ s^{- 1} \vphi < h_K \}} 
  (dd^c h_{K})^n \leq \frac{1}{s^n} \int_{\W}
 (dd^c \vphi)^n.$$
Taking an exhaustive sequence of regular compact 
subsets of the open set $ \W (s;\vphi)$ and 
using interior regularity of the capacity
we obtain our inequality in this case.\\
2)  Now for an arbitrary given function 
$\vphi \in  \mF (\W),$ there exists a decreasing 
sequence  $(\vphi_{j})$ of bounded \psh 
functions  with boundary values $0$ which
converges to $\vphi$  such that 
$\int_{\W} (dd^c \vphi)^n =
  \lim_{j} \int_{\W} (dd^c \vphi_{j})^n$ 
  (cf. [Ce 2], [Ce-Ze]).
  Then the  estimate (\ref{eq:F}) follows from the 
 first case and the lemma is proved. \\
 Now we can prove the theorem. \\
 Proof of the theorem: 1) Let $\B$ be the smallest euclidean ball of $\C^n$ containing $\W.$
Let $\vphi \in \m F (\W)$ as in the theorem and set
$\W (\vphi;s) := \{ z \in \W ; \vphi (z) \leq - s\}$ 
and $  c (s) = c_\W (s,\vphi) := \mrm{cap} (\W (\vphi;s); \W)$ for $s > 0$. Then
applying inequality 
 (\ref {eq:MAC}), we obtain
 \begin{equation}
 \la_{2 n} (\W (\vphi;s)) \leq c_n \la_{2 n} (\B) 
 \exp (- 2 c_{\W} (s)^{- 1 \slash n}),\forall s > 0.
 \label {eq:lac}
\end{equation}
Now the estimate (\ref {eq:CegC}) follows from the 
estimate (\ref{eq:lac}) using the estimate (\ref{eq:F}).

The estimate  (\ref{eq:CegR}) is proved in the 
same way using the inequalities (\ref{eq:MAR}) and
 (\ref{eq:F}).
 \fin
 \section{Global behaviour of the Lelong class}
  The next application of our theorems from the last section 
  will concern the Lelong class of \psh functions
  with logarithmic singularities at infinity defined 
  by the formula (\ref{eq:L}).

 The Lelong class of \psh functions is known  to 
 play an important role in pluripotential theory 
 (cf. [Le 1], [Be-Ta 2], [Si 1], [Si 2], [Sa 1], [Za], [Ze 1], [Ze 2]).
 
  Here we want to prove new general uniform 
 integrability theorems for the Lelong class of \psh functions.\\
 Let $g :\R^+ \to \R^+$ be an increasing function such that
 $g (0) = 0$ and $\lim_{t \to + \infty} = + \infty.$ For $\delta > 0$,
 consider the following Riemann-Stieltjes' inegral
\begin{equation}
I_\delta (g) := \delta \int_0^{+ \infty} e^{- \delta t}
 d g (t).
\label{eq:h}
\end{equation}
Then we have the following result.
 \begin{thm} 1)  For any complex euclidean closed ball 
  $\B \sub \C^n$ and any function $ u \in \mL (\C^n)$
  \begin{equation}
  \frac{1}{ \la_{2 n} (\B)} \int_{\B} g (\max_{\B} u - u) d
 \la_{2 n} \leq c_n I_2 (g),
\label{eq:h1}
\end{equation}
provided that $I_2 (g) < + \infty,$ where $c_n$ is the constant given by (\ref{eq:cn}). \\
2) Let $\G \sub \C^n$ be a generic real subspace of complex dimension $m.$
Then for any real euclidean closed ball 
  $B \sub \G$ and any function $ u \in \mL (\C^n)$
  \begin{equation}
  \frac{1}{ \la_{n + m} (B)} \int_{B} g (\max_{B} u - u) d
 \la_{n + m} \leq 8 (n + m)
  I_1(g)
\label{eq:h2}
\end{equation}
provided that $I_1 (g) < + \infty.$ 
\end{thm}
\demo We can assume $g$ to be strictly increasing. Let $\mu$ be  any Borel measure 
on $\C^n$ and $K \Sub \C^n$ any Borel set. Then for  any 
function $u \in \mL (\C^n)$ with $u\vert K \leq 0,$ we have
\begin{equation}
\int_K g (- u) d \mu = \int_0^{+ \infty} \mu(K \cap \{g (- u)\geq t\}) d t = 
\int_0^{+ \infty} \mu (K \cap \{u \leq - s\})d g (s).
\label{eq:fub}
\end{equation}
 1) Assume that $\mu := {\bf 1}_\B \la_{2 n},$ where $\B \sub
 \C^n$ is a
 complex euclidean closed ball and $u \in \mL (\C^n)$ with
 $\max_\B u = 0.$ Then by (\ref{eq:fub}), we get
 \begin{equation}
 \int_\B g (- u) d \la_{2 n}  = 
 \int_0^{+ \infty} \la_{2 n} (\B \cap \{u \leq - s\}) d g (s).
 \label{eq:fub1}
 \end{equation}
Applying the estimates (\ref{eq:lemC}) to the formula
 (\ref{eq:fub1}), we
obtain the following inequality
\begin{equation}
 \int_\B g (- u) d \la_{2 n} \leq c_n \la_{2 n} (\B)
 \int_0^{+ \infty} e^{ - 2 s} d g (s).
\label{eq:FUB2}
\end{equation}
 If  $I_2 (g) < + \infty,$ we easily see that
 $\lim_{t \to + \infty} g (t) e^{- 2 t} = 0$ and then
 by integration by parts, it follows that
 $ \int_0^{+ \infty} e^{ - 2 s} d g (s) = I_2 (g),$ which
 implies the required inequality thanks to the inequality
 (\ref{eq:FUB2}).\\
2) Assume that $\mu := {\bf 1}_B \la_{n + m},$ where
 $B \sub \G$ is a real euclidean closed ball and $u \in \mL 
(\C^n)$ with $\max_B u = 0.$ Then applying the estimates (\ref{eq:lemR}) to the formula
 (\ref{eq:fub1}), we
obtain the following inequality
 \begin{equation}
 \int_B g (- u) d \la_{n} \leq 8 (n + m) \la_{n + m} (B)
 \int_0^{+ \infty} e^{ - s} d g (s).
\label{eq:fFub3}
\end{equation}
 If $I_1 (g) < + \infty,$ then as in the first case the
 required inequality follows from the inequality
 (\ref{eq:fFub3}) by integration by parts.
\fin 
 From this general result we derive the following corollaries
 which will be useful later.
 \begin{cor}  For any complex euclidean ball 
  $\B \sub \C^n$, any function $ u \in \mL (\C^n)$ and any 
 $0 < \alpha < 2,$
  \begin{equation}
  \frac{1}{ \la_{2 n} (\B)} \int_{\B} e^{- \alpha u}
 d \la_{2 n} \leq 
  \Bigl(1 + c_n  \frac{\alpha}{2 - \alpha}\Bigr) 
  e^{- \alpha \max_{\B} u}, 
  \label{eq:ExpC}
  \end{equation}
where $c_n$ is the constant given by (\ref{eq:cn}).\\
  2) Let $\G \sub \C^n$ be a generic real subspace of
  complex dimension $m \leq n - 1.$ Then
  for any real euclidean ball 
  $B \sub \G$, any function $ u \in \mL (\C^n)$ and any $0 < \alpha < 1,$
  \begin{equation}
  \frac{1}{ \lambda_{n + m} (B)} \int_{B} e^{- \alpha u}
  d \lambda_{n + m} \leq \Bigl(1 + 8 \ (n + m) \
  \frac{\alpha}{1 - \alpha}\Bigr)
  e^{- \alpha \max_{B} u}.
  \label{eq:ExpR}
  \end{equation}
  \end {cor}
  \demo 1) Indeed,it is enough to apply Theorem 5.1 with the
 increasing function $g (t)  := e^{\alpha t} - 1,$ with $0 < \alpha < 2$ in the
 complex case and $0< \alpha < 1$ in the real
 generic case.
\fin 
\begin{cor} 1) For any complex euclidean ball 
  $\B \sub \C^n,$ any function $ u \in \mL (\C^n)$  and any
 real number 
  $p > 0,$
  \begin{equation}
  \frac{1}{\la_{2 n} (\B)} \int_{\B} (\max_{\B} u  - u)^p d
 \lambda_{2 n} \ \leq \  c_n \ 2^ {- p} \ \Gamma (p + 1),
  \label{eq:LpC}
  \end{equation}
   where $\Gamma (s) := \int_0^{+ \infty} t^{s - 1} e^{- t} d t, s > 0$ is the
   Euler function and $c_n$ is the constant given by
 (\ref{eq:cn}). \\
2) Let $\G \sub \C^n$ be a generic real subspace of
  complex dimension $m \leq n - 1.$ Then for any 
  real euclidean ball 
  $B \sub \G$, any function $ u \in \mL (\C^n)$  and any real
 number $p > 0,$
  \begin{equation}
  \frac{1}{\la_{n + m} (B)} \int_{B} 
  (\max_B u  - u)^p d 
  \la_{n + m} \ \leq \ 8 (n + m) \ \Gamma (p + 1),
  \label{eq:LpR}
  \end{equation}
 where $\Gamma$ is the Euler function.
  \end{cor}
 \demo Indeed, it is enough to apply Theorem 5.1 with the increasing function
$g (t) := t^p, t \geq 0$.
\fin 
Now we want to study the global behaviour of the
Lelong class $\mL (\C^n),$ estimating uniformly
the size of the deviation between a function and 
its mean values on complex or real euclidean balls. 

Let us recall the general definition of the space 
$\mrm{BMO}.$ Let $\G$ be a real euclidean space of 
dimension
$k \geq 1$ and let $\la_{k}$ the Lebesgue measure
on $\G.$ For a locally 
integrable function $f : \G \lgra  \ove{\R}$ and
any euclidean ball $B \sub \G,$ define the mean value
of $f$ on $B$ by
 $$f_{B} := \frac{1}{\vert B\vert_{k}} \int_{B} f 
 d \la_{k},$$
 where $\vert B\vert_{k} = \la_{k} (B).$
 Then we say that $f \in \mrm{BMO} (\G)$ if and 
 only if
 $$\Vert f\Vert_{\mrm{BMO} (\G)} := 
 \sup_{B} \{ \frac{1}{\vert B\vert_{k}} \int_{B} 
 \vert f - f_{B}\vert d \la_{k}\} < + \infty,$$ 
 where the supremum is taken over all the
 euclidean balls $B \sub \G.$ \\
   Let us first prove the following result which can be considered
   as a quatitative version for the  Lelong class 
   $\mL (\C^n)$ of 
   the classical John-Nirenberg inequality for 
   $\mrm{BMO}-$functions (cf. [St]).
\begin{thm} 1) For any complex euclidean ball 
  $\B \sub \C^n$,any function $ u \in \mL (\C^n)$  and any real number
 $\alpha < 2,$
  \begin{equation}
  \frac{1}{\vert \B \vert_{2 n}} 
  \int_{\B} e^{\alpha \vert u - u_{\B}\vert} 
  d \la_{2 n} \leq (1 + c_n \frac{\alpha}{2 - \alpha})
  \ \exp (\frac{ \alpha c_n}{2}), 
   \end{equation}
    where $u_{\B} := (1 \slash \vert \B 
    \vert_{2 n})  \int_{\B} u d \la_{2 n}$ and $c_n$ is the constant given by (\ref{eq:cn}). \\ 
  2) Let $\G \sub \C^n$ be a generic real subspace of
  complex dimension $0 \leq m \leq n - 1.$ Then for any 
  real euclidean ball  $B \sub \G$, any function $ u \in \mL (\C^n)$
  and any real number $\alpha < 1,$
   \begin{equation}
   \frac{1}{\vert B \vert_{n + m}}  
   \int_{B}  e^{\alpha \vert u - u_{B}\vert} 
   d  \la_{n + m}  \leq \Bigl(1 + 8 (n + m)
   \frac{ \alpha}{1 - \alpha}\Bigr) \ \exp ( 8 \alpha
   (n + m)), \ 
   \end{equation}
   where
 $ u_{B} := (1 \slash \vert B \vert_{n + m}) \int_{B}
  u  d \la_{n + m}.$
 \end {thm}
 \demo 1) From Corollary 5.2, it follows that for 
 a fixed function $u \in \mL (\C^n)$ and any euclidean ball
 $\B \sub \C^n,$ 
   \begin{equation}
\frac{1}{\vert \B \vert_{2 n}} 
  \int_{\B} e^{\alpha (\max_{\B} u - u)} 
  d \la_{2 n} \leq 1 +  c_ n 
  \frac{\alpha }{2 - \alpha}.
\label{eq:ExpC}
\end{equation}
  Now, from Corollary 5.3, we get
  \begin{equation}
   \max_\B u - u_\B  \leq  \frac{c_n}{2}
  \label{eq:mean-max}
  \end{equation}
   Therefore by (\ref{eq:ExpC}) and (\ref{eq:mean-max}) we get
     $$\frac{1}{\vert \B \vert_{2 n}} 
  \int_{\B} e^{\alpha \vert u - u_{\B}\vert} 
  d \la_{2 n} \leq (1 + c_n   
  \frac{\alpha }{2 - \alpha}) e^{\frac{c_n \alpha}{2}}.$$
 The real case is proved in the same way.
 \fin 
 Observe that in the complex case, a better estimate
 can be obtained using a refined version of the
 inequality (\ref{eq:mean-max}) due to Lelong 
 (cf. [Le 2], [De], [Si 2]).
 
  From the last theorem we deduce the following 
  result.
 \begin{cor} Let $\G \sub \C^n$ be a generic real subspace of
  complex dimension $m \leq n.$ Then for any
  function  $u \in \mL (\C^n)$, $ u \vert\G \in  
  \mrm{BMO} (\G)$ and 
   $$ \Vert  u \Vert_{\mrm{BMO} (\G)} 
   \leq \sigma_{n,m}, $$
In particular,  for any polynomial $P \in 
  \C [z], $ with 
 $\mrm{deg} (P) = d \geq 1,$ 
 \begin{equation}
\Vert \log \vert P \vert \Vert_{\mrm{BMO} 
  (\G)}  \leq \sigma_{n,m} \cdot d.
\label{eq:BMOC}
\end{equation}
Here $\sigma_{n,m} := 2 \log (1 + 8 (n + m)) + 8 (n + m)$  if $0 \leq m \leq n - 1$
  and   $\sigma_{n,n} := \log (1 + c_n) + c_n/2$, where $c_n$ is the constant given by (\ref{eq:cn}).
 \end{cor}
 In the totally real case where $\G = \R^n,$ the existence of
 a uniform bound for the $\mrm{BMO 
(\R^n)}-$norm of \psh functions of logarithmic
 singularities on $\C^n$ was proved earlier by A.Brudnyi 
(cf. [B 1]) with a different proof. Our proof gives a precise
quantitative estimate of the uniform bound.
\section{Global uniform integrability of \psh
 functions}
 Here we want to give a sufficient condition
 for global integrability of \psh functions  in 
 terms of the relative Monge-Amp\`ere capacity 
 of their sublevel sets.
 Then we will deduce a global integrability 
 theorem for the class of \psh functions with uniformly
 bounded
 Monge-Amp\`ere masses.

 For  any $u \in PSH^- (\W)$ and any Borel
subset $E \sub \W$ we define the truncated \psh lemniscates
associated to $u$ as 
$E (s,u) := 
\{z \in E ; u (z) <  - s\},$ for $s > 0$
and the corresponding capacity function
$$  c_{E} (s,u) =  \mrm{Cap} (E (s;u);\W).$$
 Let $\mU \sub PSH^-(\W)$ be a class of \psh 
functions on $\W$ then define
$$c_{E} (s,\mU) := \sup \{c_{E} (s,u) ; 
u \in \mU\}, s > 0.$$
 Let $g:\R^+ \to \R^+$ be an increasing function such that
 $g (0) = 0$ and $\lim_{t \to + \infty} g (t) = + \infty.$
 As in the last section, consider the following Riemann-Stieltjes'
integral for $\delta > 0,$
\begin{equation}
 I_\delta (g) := \int_0^{+ \infty}  e^{- \delta t} d g (t).
\label{eq:Ig}
\end{equation}
The main result of this section is the following.
\begin{thm} Let $\mU \sub PSH^{-}(\W)$ be a 
 class of \psh functions on $\W$ and $E \sub \W$ 
 a Borel subset 
 such that
 $$ \eta = \eta (E;\mU) := \sup_{s \geq 0} s \Bigl(
 c_{E} (s,\mU)\Bigr)^{1 \slash n} < + \infty.$$
 Then the following properties hold.\\
 1) For any function $u \in \mU$,
   $$ \int_{E} g (- u) d 
  \lambda_{2 n}  \  \leq \ c_n \tau_{2 n} (E) \ I_{2 \slash
 \eta} (g),$$
 provided that $ I_{2 \slash \eta} (g) < + \infty$
 (see \ref{eq:Ig}), where $\tau_{2 n} (E) $ is the
 $2 n-$volume
 of the smallest complex euclidean ball of $\C^n$ 
 containing $E$ and $c_n$ is the constant given by (\ref{eq:cn}). \\
 2)If $\G \sub \C^n$ is a generic real
 subspace of complex dimension $m \leq n - 1$ such that
 $ \W \cap \G \neq \es$ and $E \sub \W \cap \G$
 then for any function $ u \in \m U$,
   $$ \int_{E} g (- u) d \lambda_{n + m} 
 \ \leq \  8  (1 + \sqrt{2})
 \ (n + m) \ \tau_{n + m}  (E) \ 
 I_{1 \slash \eta} (g),$$
  provided that  $ I_{1 \slash \eta} (g) < + \infty$
 (see \ref{eq:Ig}), where  $\tau_{n + m} (E) $ is the 
 $(n + m)-$volume of the smallest euclidean ball in $\G$
 which contains $E$.  
\end{thm}
\demo 
By approximation we can assume that $g$ is strictly increasing.Let $\mu$ be any  positive Borel measure
  on $\W$ and $u \in PSH^- (\W).$ Then
  \begin{equation}
        \int_{\W} g (- u) d \mu
  = \int_{0}^{+ \infty}  \mu (\W (u;s)) 
   d g (s).
        \label{eq:Fub}
  \end{equation}
   Now let $\mu = {\bf 1}_{E} 
  \lambda_{2 n}$ and $\B$ be a complex euclidean ball of $\C^n$
containing $E.$ Then by
  (\ref{eq:MAC}) we get
 $$\la_{2 n} (E (u;s))  \leq c_n \la_{2 n} (\B) 
 \exp (- 2 c_{E} (s,u)^{- 1 \slash n}).$$
 Therefore from (\ref{eq:Fub}) we conclude that 
\begin{equation}
   \int_{E} g (- u) d \la_{2 n}
  \leq  c_n  \la_{2 n} (\B) 
   \int_{0}^{+ \infty} \exp ( - 
  2 c_{E} (s,u)^{- 1 \slash n}) d g (s).
\label{eq:la2n}
\end{equation}
  From the estimate (\ref{eq:la2n}) and the 
    hypothesis, we deduce
  that
$$ \int_{E} g (- u) d \la_{2 n}
  \leq  c_n  \la_{2 n} (\B) 
   \int_{0}^{+ \infty} \exp ( - 
  2 s \slash \eta) d g (s),
$$
which proves the required estimate.
The real generic case is proved in the same way.
\fin
From this result we can deduce the following corollaries.
\begin{cor} Let $\mU \sub PSH^{-}(\W)$ be a 
class of \psh functions on $\W$ and $E \sub \W$ be a Borel subset  such that
 $$ \eta = \eta (E;\mU) := \sup_{s \geq 0} s \Bigl(
 c_{E} (s,\mU)\Bigr)^{1 \slash n} < + \infty.$$
  Then  the following properties hold. \\
1)  For any function $u \in \mU$ and any exponent
  $0 < \alpha < 2 \slash \eta,$ 
   $$ \int_{E} e^{- \alpha u } d 
  \lambda_{2 n}  \  \leq \ \la_{2 n} (E) + c_n \tau_{2 n} (E) \ \frac{\alpha \ \eta}
  { 2 - \alpha \ \eta},$$
  where $\tau_{2 n} (E)$ is the $2 n-$volume
  of the smallest complex euclidean ball of $\C^n$ 
  containing $E$ and $c_n$ is the constant given by (\ref{eq:cn}) \\
  2) Moreover if $\G \sub \C^n$ is a generic real
  subspace of complex dimension $m \leq n - 1$
  such that $ \W \cap \G \neq \es$ and $E \sub  \W \cap \G,$ 
 for any function $ u \in \m U$ and any real number
 $\alpha < 1 \slash \eta,$ 
   $$ \int_{D} e^{- \alpha u } d  \lambda_{n + m} 
 \ \leq \ \la_{n + m} (D) +  8  (1 + \sqrt{2})
 \ (n + m) \ \tau_{n + m}  (D) \ 
  \frac{\alpha \eta}{ 1 - \alpha \eta},$$
  where  $\tau_{n + m} (D) $ is the 
  $(n + m)-$volume
  of the smallest euclidean ball of $\G$ 
  containing $D.$ 
  \end{cor} 
From the last result we can easily deduce the following one.
  \begin{cor} Let $\mU \sub PSH^{-}(\W)$ be a 
class of \psh functions on $\W$. Then the following properties hold. \\
1) If  
$$ \gamma :=
\limsup_{s \to + \infty} s \Bigl(
 c_{\W} (s,\mU)\Bigr)^{1 \slash n} < + \infty,$$
 then  for any exponent
  $0 < \alpha < 2 \slash \gamma,$ there exists a 
  constant $A_{2 n} = A_{2 n} (\alpha , \delta ,
  \W,\mU) > 0$ such that
   $$ \int_{\W} e^{- \alpha u } d 
  \lambda_{2 n}  \leq A_{2 n},  \ \forall u \in 
  \mU.$$
 2) If $\G \sub \C^n$ is a generic real subspace of
 complex dimension $m \leq n - 1$ such that
 $D := \W \cap \R^n \neq \es$ and 
$$\delta := \limsup_{s \to + \infty} s \Bigl(
 c_{D} (s,\mU)\Bigr)^{1 \slash n} < + \infty,$$ 
 then for any $\alpha < 1 \slash \delta,$ 
 there is a constant
 $A_{n,m}  = A_{n,m} (\alpha, \delta , D ,\mU) > 0$ 
 such that
   $$ \int_{D} e^{- \alpha u } d 
  \lambda_{n + m}  \leq A_{n,m} , \forall u \in \mU.$$
   \end{cor}
   \demo
  1) If $\gamma  (\mU) < + \infty,$ for any 
  $\alpha < 2 \slash \gamma (\mU),$ there is
  $s_{0} > 0$ and $\gamma_0 > 0$ such that 
  $\alpha < 2 \slash \gamma_0$ and  
 $$  s c_{\W} (s,u)^{1 \slash n} \leq \gamma_0, 
 \forall s \geq s_{0}, \forall u \in \mU.$$
  Then if we define the class $\m V := \m U + s_0,$  it
 follows that
$$  t c_{\W} (t,v)^{1 \slash n} \leq \gamma_0, 
 \forall t \geq 0, \forall v \in \m V,$$
 which implies that $\eta := \eta(\W,\m V) \leq \gamma_0.$
 Therefore, since  $\alpha < 2 \slash \gamma_0 \leq 2 \slash 
\eta,$ we can
 apply
 Theorem 6.1  to the
 class $\m V$ and get the estimate 
$$\int_{\W} e^{- \alpha v} d \lambda_{2 n} \leq
  \la_{2 n} (\W) +  c_n \tau_{2 n} (\W) \frac{\alpha \eta}{2 - \alpha \eta}.$$
 This inequality implies clearly that
$$\int_{\W} e^{- \alpha u} d \lambda_{2 n} \leq \la_{2 n} (\W) + c_n \tau_{2 n} (\W)e^{\alpha s_0}  \frac{\alpha \eta}{2 - \alpha \eta},
   \forall u \in \m U,$$
which proves the first estimate of the theorem.
The second estimate is proved in the same way.
\fin
Now we will give an application of the corollary 6.2 to the global uniform
integrability of the Cegrell class
of \psh functions of bounded Monge-Amp\`ere mass on a bounded
 hyperconvex domain.
  \begin{cor} 1) For any $\alpha < 2$ and
  any $\vphi \in \mF (\W)$ with
  $\int_{\W} (dd^c \vphi)^n \leq 1,$  
  \begin{equation}
   \int_{\W} e^{- \alpha \vphi (z)} d \la_{2n} (z) \ 
  \leq \  \la_{2n} (\W) +  c_n \ \tau_{2 n} (\W) \
  \frac {\alpha}{2 - \alpha},
  \label{eq:C-est}
  \end{equation}
 where $c_n$ is the constant given by (\ref{eq:cn}).\\
 2) If $\G \sub \C^n$ is a generic real subspace of complex
 dimension $m \leq n - 1$ such that $D := \W \cap \G \neq 
 \es,$ then for 
 any $\alpha < 1$ and any $\vphi \in \mF (\W)$
 with $\int_{\W} (dd^c \vphi)^n \leq 1,$
\begin{equation} 
\int_{D} e^{- \alpha \vphi (z)} d \la_{n + m}(z)
  \ \leq \ \la_{n + m} (D) + 8 (1 + \sqrt{2}) \
 (n + m) \
 \tau_{n + m} (D) \ \frac {\alpha}{1 - \alpha}.
  \label{eq:R-est}
  \end{equation}  
  \end{cor} 
\demo Consider the class $\mcal U$ of \psh functions in $ \vphi \mcal F (\W)$
such that $\int_{\W} (dd^c \vphi)^n \leq 1.$ Then by Lemma 4.4, we get the
inequality $\eta = \eta (E,\m U) \leq 1$ for any Borel subset $E \sub \W.$
Therefore the results above follows  immedaitely from Corollary 6.2.\fin

A uniform estimate of type (\ref{eq:C-est}) was obtained recently in ([Ce-Ze])
with a different method but with a non explicit uniform constant, while the
 estimate (\ref{eq:R-est}) seems to be new.

As in section 5, from  Theorem 6.1 we can deduce
 uniform $L^p$ estimates for functions from
 the class $\mF (\W).$
 \begin{cor} 1) For any $\vphi \in \mF (\W)$ and 
 any real number  $p > 0,$
  $$ \int_{\W} (- \vphi)^p d \lambda_{2 n} \ \leq 
 \   c_n \ \tau_{2 n} (\W) \  2^{- p} \ \Gamma (p + 1) \Bigl(\int_{\W}(dd^c
 \vphi)^n\Bigr)^{p \slash n},$$
where $c_n$ is the constant given by (\ref{eq:cn}). \\
2) If $\G \sub \C^n$ is a generic real subspace of complex
 dimension $m \leq n - 1$ such that  $D := \W \cap \G \neq 
 \es,$ then for any $\vphi \in \mF (\W)$ and any
 real number  $p > 0,$
  $$ \int_{D} (- \vphi)^p d \lambda_{n + m} \
 \leq  8 (1 + \sqrt{2})
 \ (n + m) \ \tau_{n + m} (D) \ \Gamma (p + 1) \Bigl(\int_{\W}(dd^c \vphi)^n\Bigr)^{p 
  \slash n}.$$
\end{cor}
\demo Indeed, by Lemma 4.4 the real number $\eta = \eta (E, \mcal U)$ for the
class
$\mcal U$ of plurisubharmonic functions $\vphi \in \mcal F (\W)$
such that $\int_\W (dd^c \vphi)^n \leq 1$ and any subset $E \sub \W$
satisfies the inequality $\eta \leq 1$. Since the
function $I_\delta (g)$ is decreasing in $\delta$, we easily
see that the corollary is an easy consequence of Theorem 6.1 with the function
$g (t) = t^p$.\fin 
\vskip0.3 cm
\noindent {\it Aknowlegments:} We thank Urban CEGRELL for his useful comments on the
paper and for pointing out a mistake in our earlier version of the lemma 3.3.
\vskip0.3 cm

\noindent{AMS classification:} 32U05, 32U15, 32U20, 32W20.

  \noindent Bensalem Jennane, Slimane Benelkourchi\\
 Universit\'e Mohammed V \\
 Facult\'e des Sciences de Rabat-Agdal \\
 BP. 1014, RABAT, MAROC.
\vskip 0.3 cm
\noindent Slimane Benelkourchi, Ahmed Zeriahi \\
 Universit\'e Paul Sabatier-TOULOUSE 3 \\
Institut de Math\'ematiques \\
UMR-CNRS 5580 \\
118, Route de Narbonne, F-31062 TOULOUSE. \\
e-mail: zeriahi@picard.ups-tlse.fr 

\end{document}